\documentclass[12pt]{amsart}
\usepackage{amscd,amssymb,graphicx}
\author{Alice Fialowski}
\address{E\"otv\"os Lor\'and University\\
Budapest, Hungary} \email{fialowsk@cs.elte.hu}
\author{Michael Penkava}
\address{University of Wisconsin\\
Eau Claire, WI 54702-4004} \email{penkavmr@uwec.edu}

\newtheorem{thm}{Theorem}[section]
\newtheorem{con}[thm]{Conjecture}

\theoremstyle{definition}

\def \ph{\varphi}

\def \ra{\rightarrow}

\def \hom{\mbox{\rm Hom}}
\def \ie{\hbox{\it i.e.}}

\def \tns{\otimes}
\def \mtns{\tns\cdots\tns}

\def \mcom{,\cdots,}
\def \k{\mbox{$\mathbb K$}}

\def \C{\mbox{$\mathbb C$}}
\def \Z{\mbox{$\mathbb Z$}}

\def\zt{\mbox{$\Z_2$}}

\def\inv{^{-1}}

\def\im{\operatorname{Im}}
\def\A{\mbox{$\mathcal A$}}
\def\B{\mbox{$\mathcal B$}}

\def\VA{\mbox{$\V_{\A}$}}

\def\m{\mbox{$\mathfrak m$}}

\def\V{V}

\def\dl{\delta}

\def\and{\mbox{ \rm and }}
\def\T{\mathcal T}

\def\htns{\hat\tns}

\def\htns{\hat\tns}

\def\dl#1#2{\delta^{#1}_{#2}}
\def\psa#1#2{\psi^{#1}_{#2}}

\def\TVA{\T_{\A}(\VA)}

\def\inv{^{-1}}

\def\dinfty{\mbox{$d^\infty$}}
\def\dinf{\mbox{$d^\text{inf}$}}
\def\P{\mathbb P}

\def\dl{\delta+\lambda}
\def\bdl{\bar\delta+\bar\lambda}

\def\hmu{H_\mu}

\def\ggen{\mbox{$G_{\text{gen}}$}}
\def\gdiag{\mbox{$G_\Delta$}}
\def\gdiagmd{\mbox{$\gdiag(\mu,\delta)$}}
\def\ggenmd{\mbox{$\ggen(\mu,\delta)$}}

\def\ggenmdl{\mbox{$\ggen(\mu,\delta,\lambda)$}}
\def\GL{{\mathbf{GL}}}
\newcommand{\inverselim}{\mathop{\varprojlim}\limits}

\begin{document}
\setlength{\multlinegap}{0pt}
\title{The Moduli Space of 3-dimensional Associative Algebras}%

\address{}%
\email{}%

\thanks{The research of the authors was partially supported by
 OTKA grants T043641 and T043034 and by grants from the University of Wisconsin-Eau Claire.}%

\dedicatory{Dedicated to Jim Stasheff on his 70'th birthday}%
\begin{abstract}
In this paper, we give a classification of the 3-di\-men\-sional
associa\-tive algebras over the complex numbers, including a
construction of the moduli space, using versal deformations to
determine how the space is glued together.
\end{abstract}
\maketitle

\section{Introduction}
The classification of low dimensional complex associative algebras
is one of the earliest classification theorems, dating to 1870. The
notion of deformations of associative algebras is more recent, and
was thoroughly studied by Murray Gerstenhaber in the 1960's.
Recently, the notion of moduli spaces of algebras of fixed dimension
has appeared in the literature, which gives a new, more geometric
view of the space of equivalence classes of algebras.

The classification of the algebras refers only to the set of
equivalence classes, and can be said to have been solved if
representatives for all equivalence classes have been found and
described in terms of some families and special elements. The main
concern in the classification problem is to avoid duplicate
listings.

The moduli space problem is to give a description of these
equivalence classes in terms of geometry. The set of all associative
algebra structures is described by some homogenous quadratic
polynomials on the structure constants, so these algebras determine
an algebraic set. The group of automorphisms acts on the space of
structure constants, and preserves this algebraic set. The moduli
space is the quotient by this group action, and thus is not a
variety, and as a topological space, is not Hausdorff.

In this paper, we give a different interpretation of this moduli
space, by decomposing the set of equivalence classes of associative
algebras according to deformation theory. More specifically, the
versal deformations of an element in this space give a picture of a
local neighborhood of the element. From this local analysis, the
authors were able to determine that for Lie algebras at least up to
dimension 4 (see \cite{fp8}), a geometric picture emerges. Even
though the moduli space is not Hausdorff, it has a stratification by
orbifolds, which are connected by special types of deformations,
called jump deformations, which capture all of the non-Hausdorff
behaviour of the moduli space.

Our primary aim in this paper was to examine whether the same phenomena
applied to moduli spaces of associative algebras. As the results in
this paper show, the moduli space of complex associative algebras of
dimension 3 has a stratification of the same type observed for Lie
algebras.

There are parallels and differences between the moduli spaces of
complex Lie algebras and associative algebras on a space of dimension
$n$. For example, for 2-dimensional spaces, all strata are singleton
points, but in the case of associative algebras there are six such
points, and there is only one complex Lie algebra on a space of
dimension 2.  In this paper, we shall show that on a 3-dimensional
complex space, the moduli space of complex associative algebras
consists of 22 different strata, of which 21 are singleton points, and
one an orbifold modelled on the space $\P^2/\Sigma_2$. For the
3-dimensional Lie algebra case, there were 4 strata, 3 of which were
singleton points, and one given by the orbifold $\P^2/\Sigma_2$.

\section{Early History}
Complex associative algebras of dimension up to 5 were first classified
by Benjamin Peirce as early as 1870, originally in the form of a
self-published text, which appeared later in \cite{pie}. Peirce's
methodology was criticized by his son Charles Peirce, who pointed out
that the quaternions cannot be realized as a complex associative
algebra, but are only a real algebra.  It was Charles Peirce who first
showed that the only division algebras over the reals are the real,
complex and quaternionic algebras (see \cite{pyc}). In \cite{hawkes},
an analysis of Peirce's method of construction of the complex
asssociative algebras is given, and some of the problems with the
presentation of the classification in \cite{pie} are discussed.

Since Peirce, there are other classifications of 3, 4 and 5 dimensional
associative algebras, like of B.G. Scorza, 1938 \cite{Sco}, P. Gabriel
\cite{Gab} and G. Mazzola \cite{mazz}. They all use different from ours
techniques, and we rather stayed at the original Peirce's
classification.

The main issue, which we found to be problematic when we compared our
classification with Benjamin Peirce's, is that he analyzes only those
algebras which are called by him \emph{pure algebras}, but does not
give a definition of what he means by this concept. Charles Peirce, in
footnotes added in publication, partially explains what is meant by a
non-pure algebra; more clarification is provided in \cite{hawkes}. The
fact remains that Peirce did not classify the non-pure algebras,
although we note that such a classification would not be difficult to
achieve, so Peirce's classification of the pure algebras does provide
the key insight into the classification of complex associative
algebras.  

The main technique introduced by Benjamin Peirce is his proof that
every complex algebra either contains an idempotent element, defined as
an element $a$ such that $a^2=a$, or every element is nilpotent, that
is $a^m=0$, for some $m\ge 2$. From this observation alone, he obtains
a classification of complex algebras up through dimension 5, so it is a
quite powerful tool.

Our approach to the classification problem is completely different from
the existing approaches. We use versal deformations to analyze the
moduli space, and decompose it into strata, determined by a
decomposition in terms of smooth and jump deformations. This is a new
approach to the study of moduli spaces of low dimensional algebras,
which was first introduced in \cite{fp8}.

\section{Preliminaries}
Let $V$ be a vector space over a field \k. Then an associative algebra
structure over \k\ is given by a \k-linear map $d:V\tns V\ra V$, in
other words, an element of $\hom(T^2(V),V)$. If $C^k(V)=\hom(T^k(V),V)$
is the space of $k$-cochains, then the space of cochains
$C(V)=\hom(T(V),V)$ can be expressed as a direct product
$C(V)=\prod_{k=0}^\infty C^k(V)$.

In \cite{gers}, a bracket operation, the \emph{Gerstenhaber bracket} on
$C(V)$ was introduced, which equips the space of cochains with the
structure of a \zt-graded Lie algebra, where an element of $C^k(V)$ is
odd when $k$ is even and even otherwise. Moreover $d$ is an element of
$C^2(V)$, thus odd, and $d$ is an associative algebra structure
precisely when $[d,d]=0$.

Jim Stasheff in \cite{sta4} realized that one could identify $C(V)$
with the coderivations of the tensor coalgebra $T(V)$ in a natural
manner, and that the bracket of cochains corresponds to the usual
bracket of coderivations. An odd coderivation $d$ satisfying $[d,d]=0$
is called a codifferential. Thus associative algebra structures are
simply codifferentials $d\in C^2(V)$.

Two associative algebras, given by the codifferentials $d$ and $d'$,
are isomorphic when the codifferentials are equivalent under the
following action of $G=\GL(V)$ on $C(V)$ given as follows. Any linear
map $g:V\ra V$ induces a map $g:T^k(V)\ra T^k(V)$ by $g(v_1\mtns
v_k)=g(v_1)\mtns g(v_k)$, which is bijective precisely when $g$ is
invertible. Then $g$ acts on $\ph\in C(V)$ by $g^*(\ph)=G\inv\ph g$. If
$d$ and $d'$ are the two multiplication structures, then this
definition of equivalence is the same as the usual notion of
isomorphism, that is
\begin{equation*}
d'(g(a),g(b))=g(d(a,b)).
\end{equation*}
The set of equivalence classes of codifferentials in $C^2(V)$ is called
the moduli space of associative algebra structures on $V$.

In \cite{hoch}, the notion of \emph{Hochschild cohomology} of an
associative algebra was defined. It is given by a map $D:C(V)\ra C(V)$
which maps $C^k(V)$ to $C^{k+1}(V)$. Gerstenhaber expressed this
\emph{Hochschild coboundary operator} $D$ in terms of the Gerstenhaber
bracket, by $D(\ph)=[d,\ph]$.  From this point of view, the fact that
$D^2=0$ follows immediately from the codifferential property $[m,m]=0$.

In a series of articles \cite{gers1,gers2,gers3,gers4}, Gerstenhaber
went on to use Hochschild cohomology to study the deformation theory of
associative algebras. A more general definition of deformation was
given in \cite{fi2}, in terms of a base given by a local, augmented,
commutative algebra $\A$. A local algebra is an algebra with a unique
maximal ideal $\m$, and an augmentation is a morphism
$\epsilon:\A\ra\k$. If $\A$ is a local algebra with an augmentation,
then the augmentation is unique, and is given by the map
$\epsilon\ra\A/\m\cong\k$. We shall call the algebra $\A$ infinitesimal
if $\m^2$=0, and complete if $\A=\inverselim_{n\ra\infty} \A/\m^n$. If
$V$ is a \k-vector space, then $\VA=V\tns\A$ is a free $\A$-module, the
tensor coalgebra $\TVA$ of $\VA$ over $\A$ can be identified with
$T(V)\tns\A$, and the space of cochains $C_{\A}(\VA)=\hom(\TVA,\VA)$
can be identified with $C(V)\tns\A$. (When $\A$ is a complete local
algebra, one uses $V\htns\A=\inverselim_{n\ra\infty}(V\tns\A/\m^n)$ and
$T(V)\htns\A$ instead.) Note that the augmentation induces the
structure of an $\A$-module on any \k-vector space.

A deformation of an associative algebra structure on $V$ with base $\A$
is a codifferential $d_{\A}$ on $\A$-module $V\tns\A$ , satisfying the
following property. The map $\epsilon_*:C(V)\tns\A\ra C(V)\tns\k=C(V)$,
given by $\epsilon_*(\ph\tns a)=\ph\tns\epsilon(a)$ must satisfy
$\epsilon_*(d_{\A})=d$. An automorphism of $V$ over $\A$ is an
$\A$-linear map $g_{\A}:\VA\ra\VA$, such that $\epsilon_*(g_{\A})$ is
the identity map on $V$. If $d'_{\A}=g^*_{\A}(d_{\A})$ for some
automorphism of $V$ over $\A$, then the deformation $d'_{\A}$ is said
to be equivalent to the deformation $d_{\A}$.

If $\A$ is infinitesimal, then the deformation is called an
infinitesimal deformation, and if $\A$ is complete, the deformation is
called a formal deformation. This latter terminology comes from the
fact that in the classical example of a 1-parameter formal deformation,
the deformation is given as a formal power series in a parameter.

The classical example of an infinitesimal 1-parameter deformation can
be described in terms of the infinitesimal base $A=\k[[t]]/(t^2)$, and
is given by a codifferential of the form
\begin{equation*}
d_t=d+ t\psi,
\end{equation*}
where $d$ is the codifferential giving the original multiplication. The
property that $d_t$ is a codifferential is simply $D(\psi)=0$. where
$D$ is the Hochschild coboundary operator. In terms of brackets of
cochains, this property is $[d,\psi]=0$. The classical example of a
formal 1-parameter deformation is one of the form
\begin{equation*}
d_t=d+t\psi_1+t^2\psi_2+\cdots,
\end{equation*}
with base $\k[[t]]$. The conditions which must be satisfied for a
formal deformation can be expressed in terms of brackets of the
cochains $\psi_i$, and are $\sum_{k+l=n}[\psi_k,\psi_l]=0$. for
$n=1,\dots$, where $d=\psi_0$.

If $d_{\A}$ is a deformation with base $\A$, and $f:\A\ra\B$ is a
morphism of $\k$-algebras, then there is an induced deformation
$d_{\B}$ with base $\B$, given by $d_{\B}=f_*(d_{\A})$. This allows one
to define the notion of a \emph{universal infinitesimal deformation}
$d_{\A}$, with infinitesimal base $\A$ as follows. The deformation
$d_{\A}$ is universal if, given any infinitesimal deformation $d_{\B}$
with infinitesimal base $\B$, there is a unique morphism $f:\A\ra\B$
such that $f_*(d_{\A})$ is \emph{infinitesimally equivalent} to
$d_{\B}$, in other words, equivalent as infinitesimal deformations of
$d$. In \cite{ff2}, it was proven that there is a universal
infinitesimal deformation for finite dimensional Lie algebras, over a
certain universal infinitesimal base. The proof for associative
algebras is essentially the same.

In general, there is no universal formal deformation, but there is a
weaker notion of a versal deformation. In \cite{schless} a general
framework was established for proving the existence of a versal
deformation, which was applied to Lie algebras in \cite{fi,fi2}. A
construction of a versal deformation for Lie algebras was given in
\cite{ff2}, which carries over without any difficulties for associative
algebras. A generalization of this construction to the case of infinity
algebras appeared in \cite{fp1}. A versal deformation of an associative
algebra given by the codifferential $d$ is a formal deformation
$d_{\A}$ of $d$ with base $\A$, such that if $d_{\B}$ is a formal
deformation of $d$ with base $\B$, then there is a morphism $f:\A\ra\B$
such that $f_*(d_{\A})$ is \emph{formally equivalent} to $d_{\B}$,
where formal equivalence means equivalent as formal deformations with
base $\B$. Because the morphism $f$ is not in general, unique, a versal
deformation is not universal.

The cohomology of a codifferential $d$  is $H(d)=\ker(D)/\im(D)$, where
$D$ is the Hochschild coboundary operator determined by $d$. Since
$D(C^k(V))\subseteq(C^{k+1}(V)$, we have a natural stratification of
the cohomology by spaces 
$$
H^k(d)=\ker(d:C^k(V)\ra C^{k+1}(V))/\im(d:C^{k-1}(V)\ra C^k(V)).
$$ 
Note that in deformation theory, we are only interested in calculating
cohomology with coefficients in $V$, so our notation refers to this
notion.

Let $\delta^i$ be a prebasis of $H^2(d)$. By prebasis, we mean a set of
preimages of a basis of $H^2$, lying in $C^2$, also called a set of
representative cocycles. The universal infinitesimal deformation
$\dinf$ of a codifferential $d$ is of the form
\begin{equation*}
\dinf=d+t_i\delta^i,
\end{equation*}
(we are using the summation convention on $i$), and $t_i$ are a set of
infinitesimal parameters, meaning that $t_jt_k=0$ in the infinitesimal
base $\A=\k[[t_i]]/(t_jt_k)$. Since the $\delta^i$ are cocycles, it
follows that $[\dinf,\dinf]=0$.

The infinitesimal deformation is also called the first order
deformation. To construct the second order deformation, consider the
bracket $[\dinf,\dinf]$, but this time, without assuming any property
of the parameters $t_i$.  In general, there will be second order terms
in this bracket. However, it is easy to see that $[\dinf,\dinf]$ is a
3-cocycle, so if $\alpha_i$ is a basis of $H^3$ and $\beta^i$ is a
basis of the coboundaries $D(C^2(V))$, then
\begin{equation*}
[\dinf,\dinf]=a^{ij}_kt_it_j\alpha^k+b^{ij}_kt_it_j\beta^k.
\end{equation*}
Now the $\beta_i=-\tfrac12D(\gamma_i)$ for some 2-cochains $\gamma_i$,
and if we consider the second order deformation $d^2=\dinf
+a^{ij}_kt_1t_j\gamma_k$, it is easy to see that in the bracket
$[d^2,d^2]$, the second order terms involving $\beta^k$ drop out, and
third order terms are added to the bracket.  The second order terms
involving the $\alpha^k$ do not drop out. Therefore, we assume that the
polynomials $r^k=a^{ij}_kt_it_j$ must be equal to zero, up to third
order, so for a base of the second order deformation is given by
$\k[[t_i]]/(r^k,\m^2)$, where $\m=(t_1,\dots)$ is the maximal ideal in
$\k[[t_i]]$. Continuing in this fashion, we eventually arrive at an
expression for the versal deformation
\begin{equation*}
\dinfty=d+t_i\delta^i + x_i\gamma^i
\end{equation*}
where the $x_i$ are formal power series of order at least two in the
parameters $t_i$, and the $\gamma^i$ are a prebasis of the
3-coboundaries, \ie, a preimage in $C^2$ of a basis of $B^3$. By order,
we mean the degree of the smallest term in the formal power series. We
also obtain a series of relations $r_k$, which are also given by formal
power series of order at least two, one relation for each basis element
in $H^3$. Thus the base of the formal deformations is
$\k[[t_i]]/(r_k)$. Note that the number of parameters $x_i$ is equal to
the dimension of the 3-coboundaries, while the number of parameters
$t_i$ is equal to the dimension of $H^2$.

This process of constructing the versal deformation order by order is
quite tedious, although it often terminates rather quickly. However,
the process may not terminate at all, so we have been using the
following method of computing the versal deformation.  We simply write
the formula for the versal deformation above, with undetermined
coefficients $x_i$. Consider the bracket $[\dinf,\dinf]$, and project
this to the space of 3-coboundaries. If the deformation is versal, this
projection should be equal to zero. If $M=\dim H^3(d)$, then we obtain
$M$ quadratic polynomials in the $M$ variables $x_i$, in terms of the
parameters $t_i$. In good cases, there should be a unique solution for
the $x_i$ as functions of the parameters $t_i$. Next, the projection of
$[\dinf,\dinf]$ onto $H^3$, will give $N$ polynomials in the variables
$t_i$ and $x_i$, where $N=\dim H^3$. Substituting the solutions for the
$x_i$ into these polynomials give the $N$ relations on the base of the
versal deformation. This gives the versal deformation and the relations
on the base. Note that the construction of the versal deformation
depends on the choice of prebases for $H^2$, $H^3$, $B^3$ and $B^4$ as
well as a choice of basis of $B^3$. The relations on the base and the
form of the versal deformation depend on these choices, so a lucky
choice can greatly simplify the construction.

Now, how can we use the versal deformation to construct actual
deformations, which are convergent formal power series? If one is
lucky, the relations will turn out to be rational functions of the
parameters $t_i$, and it will be possible to solve these relations to
yield simple solutions. Every such solution will determine a
deformation of $d$.  Some of the solutions are bounded away from the
origin. In this case, we say that the solution is not local, because it
is not defined for small values of the parameters. Local solutions to
the relations on the base determine how the codifferential deforms to
other codifferentials.

Our method for constructing the moduli space as a geometric object is
based on the idea that codifferentials which can be obtained by
deformations with small parameters are ``close'' to each other. From
the small deformations, we can construct 1-parameter families or even
multi-parameter families, which are defined for small values of the
parameter, except possibly when the parameters vanish. Let us consider
the 1-parameter case, since the multi-parameter case can be reduced to
studying 1-parameter families.

If $d_t$ is a one parameter family of deformations, then two things can
occur. First, it may happen that for every small value of $t$ except
zero $d_t$ is equivalent to a certain codifferential $d'$. Then we say
that $d_t$ is a jump deformation from $d$ to $d'$. It will never occur
that $d'$ is equivalent to $d$, so  there are no jump deformations from
a codifferential to itself. Otherwise, the codifferentials $d_t$ will
all be nonequivalent if $t$ is small enough. In this case, we say that
$d_t$ is a smooth deformation.

In \cite{fp10}, it was proved for Lie algebras that if one has three
codifferentials $d$, $d'$ and $d''$, and there are jump deformations
from $d$ to $d'$ and from $d'$ to $d''$, then  there is a jump
deformation from $d$ to $d''$. Similarly, if there is a jump
deformation from $d$ to $d'$, and a family of smooth deformations
$d'_t$, then there is a family $d_t$ of smooth deformations of $d$,
such that every deformation in the image of $d'_t$ lies in the image of
$d_t$, for sufficiently small values of $t$. In this case, we say that
the smooth deformation of $d$ factors through the jump deformation to
$d'$.

In the examples of complex moduli spaces of Lie algebra which we have
studied, it turns out that there is a natural stratification of the
moduli space of $n$-dimensional Lie algebras by orbifolds, where the
codifferentials on a given strata are all connected by smooth
deformations, which determine the local neighborhood structure. The
strata are connected by jump deformations, in the sense that any smooth
deformation from a codifferential on one strata to another strata
factors through a jump deformation.  Moreover, all of the strata are
given by projective orbifolds. In fact, in all the complex examples we
have studied, the orbifolds either are single points, $\P^1$,
$\P^1/\Sigma_2$ or $\P^2/\Sigma_3$. For higher dimensional complex Lie
algebras, we know that there are strata of the form $\P^n/\Sigma_{n+1}$
in the moduli space of Lie algebras of dimension $n+2$. We don't have
any proof, but we conjecture that this pattern holds in general. In
other words, we believe the following conjecture.
\begin{con}
The moduli space of Lie or associative algebras of a fixed finite
dimension $n$ are stratified by projective orbifolds, with jump
deformations and smooth deformations factoring through jump
deformations providing the only deformations between the strata.
\end{con}

\section{Extensions of Associative Algebras}

Extensions of algebraic structures have been studied and classified by
many people, for example \cite{ce,mac}. In \cite{fp9,fp11}, we gave a
treatment of this old problem that is useful in constructing
equivalence classes of extensions, because it uses cohomology as the
primary tool in the construction, so is well adapted to our
computational methods. We describe our approach to the extension
problem below.

The notion of an extension of an associative algebra $W$ by an ideal
$M$ is given by the exact sequence of algebras
\begin{equation}
0\ra M\ra V\ra W\ra 0,
\end{equation}
where $V=M\oplus W$. Let $\mu$ be the multiplication structure on $M$
and $\delta$ the multiplication structure on $W$. Then the
multiplication structure $d$ on $V$ is given by two additional
structures: the ``module structure'' $\lambda\in\hom(M\tns W\oplus
W\tns M,M)$, and the ``cocycle'' $\psi\in\hom(W\tns W,M)$. We have
$d=\delta+\mu+\lambda+\psi$, and the fact that $d$ is an associative
algebra structure is given by the codifferential property $[d,d]=0$.
This fact is equivalent to the three conditions
\begin{align*}
&[\mu,\lambda]=0\\
&[\delta,\lambda]+\tfrac12[\lambda,\lambda]+[\mu,\psi]=0\\
&[\delta+\lambda,\psi]=0.
\end{align*}
If we define $D_\mu(\ph)=[\mu,\ph]$, then $D^2_\mu=0$, so the first
condition says that $\lambda$ is a $D_\mu$-cocycle.

Under a natural notion of equivalence of extensions, cohomologous
cocycles determine equivalent extensions. This notion of equivalence is
such that if two extensions are equivalent, then the codifferentials on
the space $V$ are equivalent, so equivalent extensions determine the
same element in the moduli space of codifferentials on $V$. Denoting
the cohomology on the space of cochains $C(V)=\hom(T(V),V)$ of $V$ by
$H_\mu$, we note that the bracket on $C(V)$ descends to a bracket on
$H_\mu$. Denoting the image of a $D_\mu$-cocycle $\ph$ in $H_\mu$ by
$\bar\ph$, we note that $\bar\lambda$ and $\bar\delta$ are well
defined.

We can define an operator $D_{\delta+\lambda}$ on $H_\mu$ by
$D_{\delta+\lambda}(\bar\ph)=[\bar\delta+\bar\lambda,\bar\ph]$. By the
second condition, $D^2_{\delta+\lambda}=0$, Denote the cohomology class
of a $D_{\delta+\lambda}$-cocycle $\bar\ph$ in the cohomology
$H_{\mu,\delta+\lambda}$ induced by $D_{\delta+\lambda}$ by
$[\bar\ph]$.

We also note that $D_\mu$ commutes with $D_{\delta+\lambda}$ on $C(V)$,
which means that $D_\mu(\ker(D_{\delta+\lambda})) \subseteq
\ker(D_{\delta+\lambda})$, so induces a cohomology on this subcomplex,
which we denote by $H_\mu(\ker(D_{\delta+\lambda}))$.

Let $\gdiagmd$ be the subgroup of $\GL(M)\times\GL(W)\subseteq\GL(V)$,
consisting of those linear automorphisms $g$ of $V$ such that $g(M)=M$,
$g(W)=W$, $g^*(\mu)=\mu$ and $g^*(\delta)=\delta$. Then $g^*$ also
descends to an action on $H_{\mu}$, so  $g^*(\bar\lambda)$ is well
defined. Finally, let $\ggenmdl$ be the subgroup linear transformations
of the form $h=g\exp(\beta)$, where $g\in\gdiagmd$, and $\beta:W\ra M$,
such that $\lambda=g^*(\lambda)+[\mu,\beta]$. These are the
automorphisms of $V$ which have the property that the $\delta$, $\mu$
and $\lambda$ parts of $d$ are preserved under the action of $h^*$, and
are a subgroup of the group $\ggenmd$ of automorphisms $h$ which
preserve only $\delta$ and $\mu$.

In \cite{fp11}, the following classification of extensions was proved.
\begin{thm}
The equivalence classes of extensions of a codifferential $\delta$ on
$W$ by a codifferential $\mu$ on $M$ under the action of the group
$\ggenmd$
\begin{enumerate}
\item Isomorphism classes of $D_\mu$-cohomology classes $\bar\lambda\in H_\mu^{1,1}$ which satisfy the MC-equation
$$
\overline{[\delta+\lambda,\delta+\lambda]}=0\in  
H^{1,2}_\mu(\ker(D_{\dl}))
$$
under the action of the group $\gdiagmd$ on $\hmu$.
\item Isomorphism classes of $D_{\bdl}$-cohomology classes
$[\bar\tau]\in H_{\mu,\dl}^{0,2}$ under the action of the group
$G(\delta,\mu,\lambda)$.
\end{enumerate}
\end{thm}
Here $H_{\mu}^{1,1}$ is the subspace of $H_{\mu}$ given by the image of
elements $\lambda\in\hom(M\tns W\oplus W\tns M,M)$ and similarly,
$H_{\mu,\delta+\lambda}^{0,2}$ is the subspace given by the image of
$\tau\in\hom(W\tns W,M)$ in $H_{\mu,\delta+\lambda}$.  The cochain
$\tau$ in the theorem arises in the following manner. Given a fixed
$\lambda$, if $\psi$ is one solution to the conditions for an
extension, then the set of all such solutions is of the form
$\psi+\tau$.  The cohomology class $[\bar\tau]$ is well defined, and
the equivalence class of the extension depends only on this cohomology
class, with isomorphism classes giving equivalent extensions as stated
in the theorem.

There is a nice prescription that one can follow in applying the
theorem to construct extensions, with a minimal amount of duplication. 
Note that two extensions can be equivalent as associative algebras on
$V$ without determining equivalent extensions, so some duplication can
still occur in following the prescription.
\begin{enumerate}
\item Solve $D_\mu(\lambda)=0$. For a fixed $\mu$, this gives linear constraints on $\lambda$.
\item Solve $D_\mu(\beta)=\lambda$ for $\beta\in\hom(W,M)$. Simplify $\lambda$ by removing coboundary terms.
\item Solve the second condition on an extension. This puts linear constraints on $\psi$ and quadratic constraints
on $\lambda$.
\item Determine the group $\gdiagmd$ and apply generic elements to $\lambda$ to simplify the choice of a
representative. Now fix $\lambda$ and determine $\ggenmdl$.
\item Solve the equation $[\delta+\lambda,\psi]=0$, which gives linear constraints on $\psi$. Fix some
$\psi$ giving an extension.
\item Solve $[\mu,\tau]=[\dl,\tau]=0$ for $\tau$.
\item Determine the group $\ggenmdl$ and apply it to $\tau$ to determine the equivalence classes of $\tau$.
\end{enumerate}
To use this methodology to construct the moduli space, one simply
begins with the first choice of $\delta$ and $\mu$, determines
equivalence classes of extensions, then further checks to see which
ones are independent of the others, and then proceeds with additional
choices for $\delta$ and $\mu$.

All 3-dimensional complex associative algebras can be constructed as
extensions of a 1-dimensional algebra by a 2-dimensional ideal. This
follows from the classification of simple associative algebras.  There
are no simple 3-dimensional complex associative algebras, and thus
every 3-dimensional associative algebra must contain a 2-dimensional
ideal. Since there are 7 nonequivalent 2-dimensional associative
algebras and 2 nonequivalent 1-dimensional associative algebras, there
are 14 cases which need to be checked. The procedure is
straightforward, and we were able to check them all easily.

\section{The Moduli space of 3-dimensional Associative algebras} In
order to describe the types of associative algebras on a vector space
$W$, we introduce the following notation. Let $\{e_1\mcom e_n\}$ be a
basis of $W$. Define $\psa{ij}k:W^2\ra W$ by
\begin{equation*}
\psa{ij}k(e_me_n)=\delta^{ij}_{mn}e_k.
\end{equation*}
Then any codifferential can be expressed in the form 
$d=\psa{ij}kc_{ij}^k$, where $c_{ij}^k$ are called the structure
constants of the algebra. It is also convenient to represent $d$ in
terms of a matrix. We consider the basis of $W^2$ given by
$\{e_1^2,e_1e_2,e_2e_1,e_2^2,e_1e_3,e_2e_3,e_3e_1,e_3e_2,e_3^2,\dots\}$.
Then $d$ is given in terms of this input basis and the standard output
basis by an $n\times n^2$ matrix $A$. It is often easier to see what is
going on when studying the deformations of a codifferential to look at
the matrix of the deformed algebra.

There are only two 1-dimensional complex algebras, the trivial algebra
structure, and the structure given by the codifferential
$d_1=\psa{1,1}1$, which is the ordinary multiplication on $\C$. We will
denote the complex numbers equipped with the first structure by $\C_0$,
and with the second structure by $\C_1$. The types of 2-dimensional
complex associative algebras are given in table \ref{table 1} below.
\begin{table}[ht]\label{table 1}
\begin{equation*}
\begin{array}{lccccc}
\text{Codifferential}&H^0&H^1&H^2&H^3&H^4\\\hline\\
d_1=\psa{1,1}1&2&1&1&1&1\\
d_2=\psa{1,1}1+\psa{1,2}2&2&0&0&0&0\\
d_3=\psa{1,1}1+\psa{2,1}2&0&0&0&0&0\\
d_4=\psa{1,1}1+\psa{2,2}2&0&0&0&0&0\\
d_5=\psa{1,1}2&2&2&2&2&2\\
d_6=\psa{1,1}1+\psa{1,2}2+\psa{2,1}2&2&1&1&1&1\\\\\hline
\end{array}
\end{equation*}
\caption{\protect Cohomology of two dimensional complex associative
algebras}\label{Table 1}
\end{table}

 The classification of the algebras can already be
determined from \cite{pie}, which explicitly gives $d_2$, $d_5$ and
$d_6$, and implicitly $d_3$. The algebras $d_1$ and $d_4$ are
decomposable as direct sums.  We will denote $\C^2$, equipped with the
algebra structure given by the codifferential $d_i$ as $\C^2_i$. Thus
$\C^2_1=\C_1\oplus \C_0$, and $\C^2_4=\C_1\oplus\C_1$ gives a
representation of the structures $d_1$ and $d_4$ as direct sums.

Of these algebras, only $d_2$ and $d_3$ are not commutative.  In terms
of Peirce's classification by nilpotent and idempotent elements, types
$d_2$, $d_3$ and $d_5$ contain nonzero idempotents, while $d_5$ is a
nilpotent algebra. The algebras $d_4$ and $d_6\cong \C[x]/(x^2)$ are
the only unital algebras.

\subsection{Three dimensional algebras}
A complete list of the 22 distinct types of 3-dimensional complex
associative algebras is given in Table \ref{Table 2}.
\begin{table}[ht]
\begin{equation*}
\begin{array}{lrrrrr}
\text{Codifferential}&H^0&H^1&H^2&H^3&H^4\\\hline\\
d_1=\psa{1,1}1&3&4&8&16&32\\
d_2=\psa{1,1}1+\psa{2,2}3&3&2&2&2&2\\
d_3=\psa{1,1}1+\psa{1,3}3&1&1&2&2&2\\
d_4=\psa{1,1}1+\psa{3,1}3&1&1&2&2&2\\
d_5=\psa{1,1}1+\psa{1,3}3+\psa{3,1}3&3&2&2&2&2\\
d_6=\psa{1,1}1+\psa{3,3}3&3&1&1&1&1\\
d_7=\psa{1,1}1+\psa{2,1}2+\psa{1,3}3&0&1&0&1&0\\
d_8=\psa{1,1}1+\psa{2,1}2+\psa{3,1}3&0&3&0&0&0\\
d_9=\psa{1,1}1+\psa{2,1}2+\psa{1,3}3+\psa{3,1}3&1&1&1&1&1\\
d_{10}=\psa{1,1}1+\psa{2,1}2+\psa{3,3}3&1&0&0&0&0\\
d_{11}=\psa{1,1}1+\psa{2,2}2+\psa{2,3}3&1&0&0&0&0\\
d_{12}=\psa{1,1}1+\psa{2.2}2+\psa{2,3}3+\psa{3.2}3&3&1&1&1&1\\
d_{13}=\psa{1.1}1+\psa{2,2}2+\psa{2,3}3+\psa{3,1}3&1&0&0&0&0\\
d_{14}=\psa{1,1}1+\psa{2,2}2+\psa{3,3}3&3&0&0&0&0\\
d_{15}=\psa{1,1}2&3&5&9&17&33\\
d_{16}=\psa{1.1}2+\psa{1,2}3+\psa{2,1}3&3&3&3&3&3\\
d_{17}=\psa{1,1}1+\psa{1,1}2+\psa{1,2}2+\psa{2,1}2+\psa{1,3}3&1&1&1&1&1\\
d_{18}=\psa{1,1}1+\psa{1,1}2+\psa{1,2}2+\psa{2,1}2+\psa{1,3}3+\psa{3,1}3&3&4&6&12&24\\
d_{19}=\psa{3,3}3+\psa{1,1}2+\psa{1,3}1+\psa{3,1}1+\psa{2,3}2+\psa{3,2}2&3&2&2&2&2\\
d_{20}=\psa{1,1}1+\psa{1,2}2+\psa{1,3}3&3&0&0&0&0\\
d_{21}=\psa{1,1}3+\psa{1,2}3-\psa{2,1}3&1&2&3&4&5\\
d_{22}(1:0)=\psa{1,2}3&1&2&5&8&11\\
d_{22}(1:1)=\psa{1,2}3+\psa{2,1}3&3&4&5&7&8\\
d_{22}(1:-1)=\psa{1,2}3-\psa{2,1}3&1&4&5&8&9\\
d_{22}(x:y)=x\psa{1,2}3+y\psa{2,1}3&1&2&2&3&4\\\hline
\end{array}
\end{equation*}
\caption{\protect Three dimensional complex algebras and their cohomology}\label{Table 2}
\end{table}
There is a family of algebras labeled $d_{22}(x:y)$, where the
parameter $(x:y)$ in $d_{22}$ is a projective coordinate,
representing the fact that any nonzero multiple of the parameter
gives an equivalent deformation. Actually, $d_{22}(x:y)\sim
d_{22}(y:x)$, so that this family of differentials is parameterized
by $\P^1/\Sigma_2$, where the action of the permutation group
$\Sigma_2$ on $\P^1$ is by interchanging coordinates.

The codifferential $d_{19}$ is expressed as an extension of a
$2$-dimensional algebra by a $1$-dimensional algebra.
However, it is equivalent to the codifferential
\begin{equation*}
d=\psa{1,1}1+\psa{1,1}2+\psa{1,2}2+\psa{2,1}2+\psa{1,2}3+\psa{2,1}3+\psa{2,2}3+\psa{1,3}3+\psa{3,1}3,
\end{equation*}
which is an extension of a $1$-dimensional algebra by a
$2$-di\-men\-sional algebra. The reason it is given in the form in
which it appears in the table is that the cohomology was much
simpler to compute in this form. However, the fact that it is
equivalent to the codifferential above implies that every
3-dimensional associative algebra has a 2-dimensional ideal.

A partial classification of 3-dimensional associative algebras
appears already in \cite{pie}. Our method of obtaining these
codifferentials differs from Peirce's in that we constructed them by
studying extensions of $2$-dimensional associative algebras by
$1$-dimensional algebras (and vice-versa). By Wedderburn's theorem,
there are no simple 3-dimensional complex associative algebras, so
all three dimensional complex associative algebras arise as
extensions.

The algebras $d_1\mcom d_6$, as well as $d_{10}\mcom d_{13}$,
$d_{14}$ and $d_{15}$ are direct sums, while the algebras $d_7\mcom
d_9$, $d_{13}$, $d_{17}$, $d_{18}$ and $d_{20}$ are examples of what
Peirce calls not pure algebras. The algebras $d_{16}$, $d_{19}$,
$d_{21}$ and $d_{22}$ all correspond to algebras on Peirce's list.

In a footnote in \cite{pie}, added by Peirce's son Charles Peirce,
there is an explanation that a not pure algebra is one in which
there is a basis $a$, $b$ and $c$ such that $a$ and $b$ span a
subalgebra, $a$ and $c$ span another, and $bc=cb=0$. If one assumes
that products of $a$ and $b$ are multiples of $b$ and that products
of $a$ and $c$ are multiples of $c$, then one obtains automatically
the structure of an associative algebra, whenever the two
subalgebras are associative. Of course, the algebras which are
direct sums are also not pure in this sense.

Peirce does not give a classification of the algebras which are not
pure, confining himself to classifying the pure algebras. However,
it would probably not be a difficult task to classify the not pure
algebras using the methods in Peirce's paper. Thus, our analysis
concurs with Peirce's classification.

To compute the nonequivalent codifferentials, we used the methods
outlined in \cite{fp11}, which outlines a method of computing
nonequivalent extensions using cohomological methods.
\section{Versal Deformations of the codifferentials}
\subsection{Type $d_1$}
The codifferential $d_1=\psa{1,1}1$ represents the multiplication structure $\C_1\oplus\C_0\oplus\C_0$.
Since $\dim H^2=8$ and $\dim H^3=16$, one might expect that the versal deformation would be quite complicated.
However, it turns out that the infinitesimal deformation coincides with the versal deformation. The matrix of the
versal deformation is given by
\begin{equation*}
 \left[ \begin {array}{ccccccccc} 1&0&0&0&0&0&0&0&0
\\\noalign{\medskip}0&0&0&t_{{8}}&0&t_{{6}}&0&t_{{4}}&t_{{3}}
\\\noalign{\medskip}0&0&0&t_{{7}}&0&t_{{5}}&0&t_{{2}}&t_{{1}}
\end {array} \right].
\end{equation*}
There are 16 relations on the base of the versal deformation, which we omit, for sake of brevity.
Although the relations are numerous, their solution is not so complicated. Using Maple, we found the following seven
solutions to the relations.
\begin{align*}
& t_3=t_5=t_4=t_7=0,t_{{8}}=t_{{2}},t_{{1}}=t_{{6}}\\
& t_3=t_4=t_6=0,t_{{1}}=-{\frac {t_{{5}} \left( -t_{{5}}+t_{{8}} \right) }{t_{{7}}}},t_{{2}}=t_{{5}} \\
&t_3=t_5=t_2=0, t_{{1}}=t_{{6}},t_{{4}}=t_{{6}}\\
& t_2=t_3=t_6=t_7=0,t_{{8}}=t_{{5}},t_{{1}}=t_{{4}}\\
& t_2=t_5=t_7=t_8=0,t_{{1}}=t_{{6}},t_{{4}}=t_{{6}}\\
&t_2=t_5=t_7=0, t_{{6}}=t_{{4}},t_{{3}}=-{\frac {t_{{4}} \left( t_{{1}}-t_{{4}} \right) }{t_{{8}}}}\\
&t_{{2}}=t_{{5}},t_{{6}}=t_{{4}},t_{{8}}={\frac {-t_{{1}}t_{{4}}+t_{{3}}t_{{5}}+{t_{
{4}}}^{2}}{t_{{3}}}},t_{{7}}={\frac {t_{{4}}t_{{5}}}{t_{{3}}}}
\end{align*}
In order to determine which codifferential a deformation is equivalent
to, one needs to substitute a solution to the relations into the
general form for the versal deformation. Thus we studied 7 versions of
the versal deformation, corresponding to the 7 solutions above. We omit
giving details about which codifferentials each of the 7 turns out to
be equivalent to. What is important is which codifferentials arise as
deformations for small values of the parameters. Note each of the
solutions above remains valid if we substitute all of the parameters
with zero. Such a solution will be called \emph{local}. Only the local
solutions are relevant to the description of how the moduli space is
glued together.  We found that $d_1$ deforms to types $d_2$, $d_6$,
$d_{10}$, $d_{11}$, $d_{12}$, and $d_{14}$.

\subsection{Type $d_2$}
The codifferential $d_2$ represents the multiplication structure $\C_1\oplus\C^2_5$. Since $\dim H^2=\dim H^3=2$,
the versal deformation is not so difficult to compute, and it turns out that the infinitesimal deformation
is versal in this case as well; moreover, the relations on the base vanish, so that the coderivation
$\dinfty$, represented by the matrix
\begin{equation*}
\left[ \begin {array}{ccccccccc} 1&0&0&0&0&0&0&0&0\\\noalign{\medskip}0&0&0&t_{{1}}&0&t_{{2}}&0&t_{{2}}&0
\\\noalign{\medskip}0&0&0&1&0&0&0&0&t_{{2}}\end {array} \right]
\end{equation*}
is indeed a codifferential for all values of the parameters. We found that $d_2$ deforms to $d_6$, $d_{12}$
and $d_{14}$.

\subsection{Type $d_3$}
The codifferential represents the multiplication structure $\C^2_2\oplus\C_0$. We have $\dim H^2=\dim H^3=2$,
and the versal deformation is given by the infinitesimal deformation. The matrix of $\dinfty$ is given by
\begin{equation*}
\left[ \begin {array}{ccccccccc} 1&0&0&0&0&0&0&0&0\\\noalign{\medskip}0&0&0&t_{{2}}&0&0&0&0&0\\\noalign{\medskip}0&0&0&0
&1&0&0&t_{{1}}&0\end {array} \right].
\end{equation*}
There is one nontrivial relation $t_1(t_2-t_1)=0$, which has the two solutions $t_1=0$ and $t_2=t_1$, both
of which are local. We found that $d_3$ deforms to $d_{11}$ and $d_{13}$.
\subsection{Type $d_4$}
The codifferential $d_4$ represents the multiplication structure $\C^2_3\oplus\C_0$. The infinitesimal deformation
is versal, and is given by the matrix
\begin{equation*}
\left[ \begin {array}{ccccccccc} 1&0&0&0&0&0&0&0&0\\\noalign{\medskip}0&0&0&t_{{1}}&0&0&0&0&0\\\noalign{\medskip}0&0&0&0
&0&t_{{2}}&1&0&0\end {array} \right].
\end{equation*}
There is one nontrivial relation $t_2(t_1-t_2)=0$. We found that $d_4$ deforms to types $d_{10}$ and $d_{13}$.
Note that $d_3$ and $d_4$ are mirror images of each other, so they should have the same type of pattern to
their deformations.
\subsection{Type $d_5$}
The codifferential $d_5$ represents the multiplication structure $C^2_6\oplus C_0$. The infinitesimal deformation,
given by the matrix
\begin{equation*}
\left[ \begin {array}{ccccccccc} 1&0&0&0&0&0&0&0&t_{{1}}\\\noalign{\medskip}0&0&0&t_{{2}}&0&0&0&0&0\\\noalign{\medskip}0&0&0&0
&1&0&1&0&0\end {array} \right],
\end{equation*}
is versal and there are no relations.  We found that $d_5$ deforms to types $d_6$, $d_{12}$ and $d_{14}$.
\subsection{Type $d_6$}
The codifferential $d_6$ represents the multiplication structure $\C^2_4\oplus\C_0=\C_1\oplus \C_1\oplus \C_0$. The matrix of
the infinitesimal deformation, which is versal, is given by
\begin{equation*}
\left[ \begin {array}{ccccccccc}
1&0&0&0&0&0&0&0&0\\\noalign{\medskip}0&0&0&t_1&0&0&0&0&0\\\noalign{\medskip}0&0&0&0&0&0&0
&0&1\end {array} \right].
\end{equation*}
There are no relations on the base.  We found that $d_6$ deforms only to type $d_{14}$.
\subsection{Type $d_7$}
The codifferential $d_7$ does not decompose as a direct sum. Since $H^2=0$, there are no deformations.
Note that the cohomology vanishes in even degrees, and has dimension $1$ in odd degrees. Therefore, there
are no odd elements.
\subsection{Type $d_8$}
The codifferential $d_8$ does not decompose as a direct sum. Since $H^2=0$, there are no deformations.
In this case, only $H^1$ does not vanish.
\subsection{Type $d_9$}
The codifferential $d_9$ does not decompose as a direct sum. Its infinitesimal deformation, which is versal,
is given by the matrix
\begin{equation*}
\left[ \begin {array}{ccccccccc} 1&0&0&0&0&0&0&0&0\\\noalign{\medskip}0&0&1&0&0&0&0&0&0\\\noalign{\medskip}0&0&0&0&1&0&1
&0&t_{{1}}\end {array} \right].
\end{equation*}
We found that $d_9$ deforms only  to $d_{10}$.
\subsection{Type $d_{10}$}
The codifferential $d_{10}$ represents the multiplication structure $\C^2_3\oplus\C_1$. Since $H^2=0$, it
does not deform.
\subsection{Type $d_{11}$}
The codifferential $d_{11}$ represents the multiplication structure $\C^2_2\oplus\C_1$. Since $H^2=0$ it does not
deform.
\subsection{Type $d_{12}$}
The codifferential $d_{12}$ does not decompose as a direct sum. Its infinitesimal deformation, given by
the matrix
\begin{equation*}
\left[ \begin {array}{ccccccccc} 1&0&0&0&0&0&0&0&0\\\noalign{\medskip}0&0&0&1&0&0&0&0&t_{{1}}\\\noalign{\medskip}0&0&0&0
&0&1&0&1&0\end {array} \right],
\end{equation*}
is versal, and there are no relations.  We found that $d_{12}$ deforms only to $d_{14}$.
\subsection{Type $d_{13}$}
The codifferential $d_{13}$ does not decompose as a direct sum, and it has no deformations, since $H^2=0$.
\subsection{Type $d_{14}$}
The codifferential $d_{14}$ decomposes as $\C^2_4\oplus\C_1=\C_1\oplus\C_1\oplus\C_1$. There are no deformations,
since $H^2=0$.
\subsection{Type $d_{15}$}
The codifferential $d_{15}$ represents the multiplication structure $\C^2_5\oplus\C_0$. Since the dimension
of $H^2$ is 9, the largest dimension of $H^2$ for any of the codifferentials, it is not surprising that
the versal deformation is quite complex. Its infinitesimal deformation is not versal.
We shall omit giving the matrix of the versal deformation, but mention that it is given in terms of at most
cubic polynomials in the parameters, so that the third order deformation is versal.
We also omit the 17 relations on the base of the versal deformation. We note that with the help of Maple, we
decomposed the relations into 9 solutions, all of which are local.

We found that
$d_{15}$ deforms to types $d_1$, $d_2$, $d_3$, $d_4$, $d_5$, $d_6$,
$d_9$, $d_{10}$, $d_{11}$, $d_{12}$, $d_{13}$, $d_{14}$,
$d_{16}$,$d_{17}$, $d_{18}$, $d_{19}$,
$d_{21}$ and $d_{22}(x:y)$ for every value of the parameter $(x:y)$ except $(1:-1)$.
In other words, it deforms to all of
the codifferentials except for $d_7$, $d_8$,  $d_{20}$, and $d_{22}(1:-1)$.
\subsection{Type $d_{16}$}
The codifferential $d_{16}$ is not a direct sum, and one of the
algebra's given in Peirce's article. The infiniitesimal deformation is versal and has matrix
\begin{equation*} \left[ \begin {array}{ccccccccc} -t_{{2}}&t_{{3}}&t_{{3}}&t_{{1}}&t_{
{1}}&0&t_{{1}}&0&0\\\noalign{\medskip}1&0&0&t_{{3}}&0&t_{{1}}&0&t_{{1}
}&0\\\noalign{\medskip}0&1&1&t_{{2}}&0&0&0&0&t_{{1}}\end {array}
 \right]
\end{equation*}
We found
that $d_{16}$ deforms to $d_2$, $d_5$, $d_6$, $d_{12}$ $d_{14}$ and $d_{19}$.

\subsection{Type $d_{17}$}
The codifferential $d_{17}$ is not a direct sum.
We found that $d_{17}$ deforms only to type $d_{11}$.
\subsection{Type $d_{18}$}
The codifferential $d_{18}$ is not a direct sum.  Since $\dim
H^2=6$  and $\dim H^3=12$, it is not surprising that there might
be some difficulty computing the versal deformation, and indeed,
we were unable to compute the relations on the base of the versal
deformation explicitly. Nevertheless, we were able to determine
all the local solutions to the relations and determine which
codifferentials they were equivalent to. We found that $d_{18}$
deforms to types $d_{12}$, $d_{13}$, $d_{14}$ and $d_{19}$.

\subsection{Type $d_{19}$}
The codifferential $d_{19}$ is not a direct sum, and it appears in Peirce's classification. Its infinitesimal
deformation is versal, and is given by the matrix
\begin{equation*}
 \left[ \begin {array}{ccccccccc} 0&0&0&t_{{2}}&1&0&1&0&0
\\\noalign{\medskip}1&0&0&t_{{1}}&0&1&0&1&0\\\noalign{\medskip}-t_{{1}
}&t_{{2}}&t_{{2}}&0&0&0&0&0&1\end {array} \right].
\end{equation*}
We found that $d_{19}$ deforms  to types $d_{12}$ and $d_{14}$.
\subsection{Type $d_{20}$}
The codifferential $d_{20}$ is not a direct sum. Since $H^2=0$, it has no deformations.
\subsection{Type $d_{21}$}
The codifferential $d_{21}$ is not a direct sum, and it appears in Peirce's classification.
The relations on the base of the versal deformation are
\begin{align*}
&{\frac { \left( t_{{3}}+1 \right)  \left( -t_{{1}}+t_{{2}} \right)
 \left( 2\,t_{{1}}-t_{{2}}+t_{{2}}t_{{3}} \right) t_{{3}}}{
 \left( -3+t_{{3}} \right)  \left( t_{{3}}-1 \right) }}=0\\
& {\frac {
 \left( 2\,t_{{1}}-t_{{2}}+t_{{2}}t_{{3}} \right)  \left( 2\,t_{{1}}t_
{{3}}-2\,t_{{1}}-t_{{2}}-4\,t_{{2}}t_{{3}}+t_{{2}}{t_{{3}}}^{2}
 \right) }{ \left( -3+t_{{3}} \right)  \left( t_{{3}}-1 \right) }}=0\\
&\,{\frac { \left( -t_{{1}}+t_{{2}} \right)  \left( 2\,t_{{1}}-t_{{2}}+
t_{{2}}t_{{3}} \right) }{-3+t_{{3}}}}=0\\
&{\frac { \left( -t_{{1}}+t_{{2}
} \right)  \left( t_{{1}}+t_{{2}} \right)  \left( t_{{3}}+1 \right) }{
-3+t_{{3}}}}=0.
\end{align*}
Note that these expressions are rational, not polynomial, in the $t$-parameters. As a consequence, calculation of
the versal deformation order by order would not be successful. However, there are only two solutions to these
equations, given by
\begin{align*}
&t_3=-1, t_2=t_1\\
&t_1=t_2=0\\
\end{align*}
The first of these is not local. The other solution give a very simple matrix for the versal deformation:
\begin{equation*}
\left[ \begin {array}{ccccccccc} 0&0&0&0&0&0&0&0&0\\\noalign{\medskip}0&0&0&0&0&0&0&0&0\\\noalign{\medskip}1&1&-1&t_{{3}
}&0&0&0&0&0\end {array} \right].
\end{equation*}
Thus, even though $H^2$ has dimension 3, effectively it as is if the dimension was only 1, because the relations
force such a simplification in any true deformation.

We found that $d_{21}$ deforms to $d_{22}(x:y)$, except for the values $(1:0)$, $(1:1)$ and $(1:-1)$. The exceptions
turn out to be important, because for these exceptional values of the parameter $(x:y)$, the codifferential
$d_{22}$ has extra deformations.  Since $d_{21}$ does not have any extra deformations, it would be impossible
for it to deform into these types.  This is not the reason we know it has no such deformations; rather we
were gratified to see that our computations of the deformations of $d_{21}$ did not lead to any contradictions.
\subsection{Type $d_{22}(x:y)$}
The codifferential $d_{22}(x:y)$ is not a direct sum, and appears as a family on Peirce's list. More precisely,
Peirce gives a 1-parameter family of associative algebra structures, parameterized by $a$, which mostly coincides
with our family.  The case $a=-2$, which Peirce points out as being special, coincides with our codifferential
$d_{19}$. Also $a=2$ corresponds to the codifferential $d_{19}$, and more generally, the multiplication determined
by any parameter $a$ is equivalent to that determined by $-a$, although this fact was not mentioned in
\cite{pie}.  The coderivations $d_{22}(1:0)$ and $d_{22}(1:-1)$ are not equivalent to members of Peirce's family,
but all the rest are.

As was mentioned in the beginning of the section, the family
$d_{22}(x:y)$ is parameterized by $(x:y)\in\P^1/\Sigma_2$. Note that
Peirce's family is parameterized by $\C/\Sigma_2$, where $\Sigma_2$
acts by taking $a$ to $-a$. Topologically, the resulting quotient is
still $\C$, and if we remove the point corresponding to $a=2$, we
obtain $\P^1$ minus two points. Since the orbifold $\P^1/\Sigma_2$
which parameterizes $d_{22}$ is topologically $\P^1$, removing the
points $(1:0)$ and $(1:-1)$ gives the same topological type as
Peirce's family.

Generically, the dimension of $H^3$ is 3, but for the special cases $(1:0)$, $(1:1)$, and $(1:-1)$ the cohomology
has dimension 5, so  these cases require separate treatment. Let us describe the generic case first. The
matrix of the versal deformation is complicated, so we omit it. The relations on the base of the versal
deformation are
\begin{align*}
&{\frac { \left( r+t_{{1}} \right)  \left( t_{{1}}+s+r \right)
t_{{3}}t_{{2}}}{s \left( r+t_{{1}}-s \right) }}=0\\
&{\frac { \left( r+t_{{1}} \right) ^{2} \left( t_{{1}}+s+r \right) {t_{{3}}}^{2}}{{s}^{2}
 \left( r+t_{{1}}-s \right) }}=0\\
&{\frac { \left( t_{{1}}+s+r
 \right) {t_{{2}}}^{2}{s}^{2}}{ \left( r+t_{{1}} \right) ^{2} \left( r
+t_{{1}}-s \right) }}=0,
\end{align*}
with two solutions:
\begin{align*}
t_2=t_3=0\\
t_1=-(x+y).
\end{align*}
The second solution is not local, unless $y=-x$, which is the special case $d(1:-1)$ to be considered later.
Thus, we can restrict our attention to the first solution, in which case  the versal deformation is simply
$\dinfty=d_{22}(x+t_1:y)$.  This is the first case we have studied in which a deformation occurs which is not
a jump deformation.  This smooth deformation is simply a deformation along the family. Thus generically,
a codifferential $d_{22}(x:y)$ only deforms along the family $d_{22}$, and it deforms smoothly along the
family.
\subsubsection{The codifferential $d_{22}(1:0)$}
For this special case, the base of the versal deformation is given by 8 relations, two of which are trivial.
Solving the relations, one obtains 6 solutions, only 4 of which turn out to be local.
We found that $d_{22}(1:0)$ deforms to types $d_3$, $d_4$, $d_9$, $d_{10}$, $d_{11}$, $d_{13}$ and $d_{17}$,
as well as smoothly deforming along the family.
\subsubsection{The codifferential $d_{22}(1:1)$}
For this special case, the base of the versal deformation is given by 7 relations, which yield 3 solutions,
2 of which are local. We found that $d_{22}(1:1)$ deforms to types
$d_2$, $d_5$, $d_6$, $d_{12}$, $d_{14}$, $d_{16}$ and $d_{19}$. Once you know that it deforms to type
$d_{16}$, this determines all the rest, because $d_{16}$ deforms to all the other ones.  We also found that
$d_{22}(1:1)$ deforms along the family $d_{22}$.
\subsubsection{The codifferential $d_{22}(1:-1)$}
For this special case, the base of the versal deformations is given by 8 relations, which are very complicated,
but still only give 3 solutions, 2 of which are local. As well as deforming along the family, we find that
$d(1:-1)$ also deforms to type $d_7$.
\section{Gluing the moduli space together}

Let us define the level of an algebra as follows. A codifferential is on level 0 if it has no deformations.
A codifferential is on level $k+1$ if some element of its family
has a deformation to an element on level $k$, and other than deformations
along its own family, all deformations are to elements of level $k$ or below. Using this definition of level, there
are 7 levels for the elements in our moduli space.

On level 0, there are 7 elements: $d_7$, $d_8$, $d_{10}$, $d_{11}$, $d_{13}$, $d_{14}$ and $d_{20}$. Of these,
$d_8$ and $d_{20}$ don't even have any jump deformations to them from higher level objects, so they are truly
isolated.

On level 1, there are 6 elements: $d_3$, $d_4$, $d_6$, $d_9$,
$d_{12}$ and $d_{17}$. We have $d_3$ deforming to both $d_{11}$
and $d_{13}$, $d_4$ to both $d_{10}$ and $d_{13}$, $d_6$ deforms
to $d_{14}$ only, $d_9$ deforms only to $d_{10}$, $d_{12}$ to
$d_{14}$ only, and $d_{17}$ deforms only to $d_{11}$.

On level 2, there are 3 elements: $d_2$, $d_5$ and $d_{19}$. Both
$d_2$ and $d_5$ deform to $d_6$ and $d_{12}$ on level 1, and
$d_{19}$ deforms to $d_{12}$ on level 1.

On level 3, there are 3 elements: $d_1$, $d_{16}$ and $d_{18}$.  Here the picture is more complex, because the
jumps to lower levels do not always factor through a jump to level 2.  We have $d_1$ deforming to $d_2$ on level 2,
as well as directly to $d_{10}$ and $d_{11}$ on level 0. The behaviour of $d_{16}$ is determined by jumps to
$d_2$, $d_5$ and $d_{19}$ on level 2. Finally, $d_{18}$ deforms to $d_{19}$ on level 2, as well as to $d_{13}$ on
level 0.

On level 4, there is only the family $d_{22}(x:y)$. While most elements deform only along the family, the special
element $d_{22}(1:1)$ deforms to $d_{16}$ on level 3, which puts the family on level 4. In addition, $d_{22}(1:0)$
deforms to $d_3$, $d_4$, $d_9$ and $d_{17}$ on level 1, while $d_{22}(1:-1)$ has a jump to $d_7$ on level 0.

On level 5, there is only the element $d_{21}$, which has jump deformations to most elements of the family $d_{22}(x:y)$
on level 4.

Finally, on level 6, there is only $d_{15}$, which has a jump to $d_{21}$ on level 5, as well as jumps to almost all
the other codifferentials.

The moduli space of three dimensional complex associative algebras is illustrated in Figure \ref{Figure1} below .

\begin{figure}[ht]
\vspace{.05in}
\includegraphics[width=4.5in]{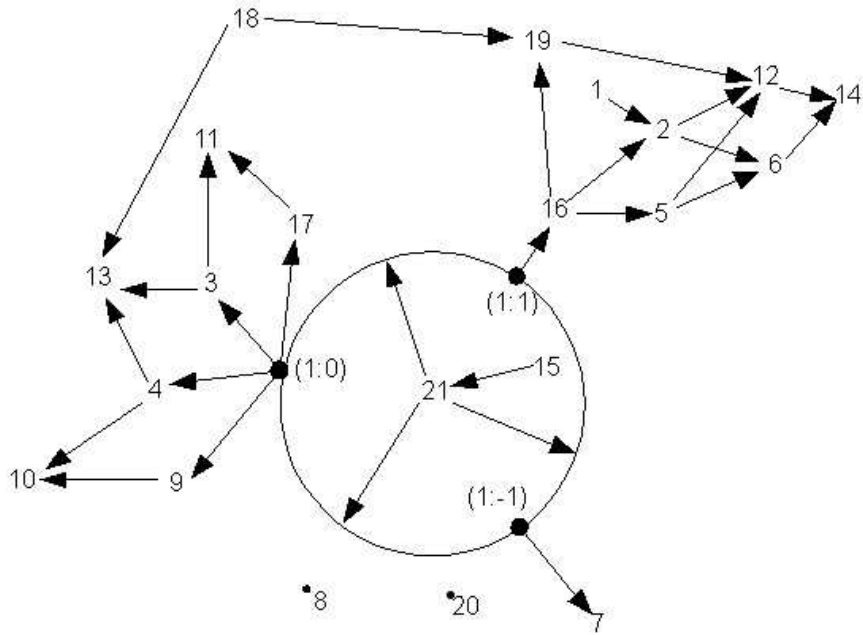}

\caption{The Moduli Space of 3-dimensional associative algebras. } \label{Figure1}
\end{figure}

\section{Conclusions} The authors have been studying moduli spaces of low dimensional Lie algebras, and have noticed that
the description of the moduli space can be given in terms of a
stratification by complex projective orbifolds, with families
given by very special orbifolds of the type $\P^n/\Sigma_{n+1}$.
There are a lot of similarities between the picture for three
dimensional associative algebras and the Lie algebra case. There
is one family, parameterized by $\P^1/\Sigma_2$, and some special
points, for both the three dimensional Lie and the three
dimensional associative algebra cases. However, there are only 6
special points in the Lie case, and 21 special points in the
associative algebra case. In the Lie case, there are no elements
which have the property that they neither deform, nor does any
element deform to them, but in the associative case, there are two
such elements. So there are similarities and differences in the
pictures of Lie and associative algebras.

\section{Acknowledgments} The authors would like to thank Jim Stasheff,
who read this manuscript and made many helpful remarks.  In addition,
the authors thank Mitch Phillipson, an undergraduate student, for
creating the illustration in Figure 1. We also thank the referee for
providing useful historical references. 

\bibliographystyle{amsplain}

\providecommand{\bysame}{\leavevmode\hbox to3em{\hrulefill}\thinspace}
\providecommand{\MR}{\relax\ifhmode\unskip\space\fi MR }
\providecommand{\MRhref}[2]{%
  \href{http://www.ams.org/mathscinet-getitem?mr=#1}{#2}
}
\providecommand{\href}[2]{#2}

\end{document}